\documentclass[12pt]{article}
\usepackage{amsmath}
\usepackage{amssymb}
\allowdisplaybreaks

\def\q{\quad}

\def\qtq#1{\q\t{#1}\q}
\def\mod#1{\ (\text{\rm mod}\ #1)}
\def\t{\text}
\def\f{\frac}
\def\e{\equiv}
\def\b{\binom}
\def\sls#1#2{(\f{#1}{#2})}
 \def\ls#1#2{\big(\f{#1}{#2}\big)}
\def\Ls#1#2{\Big(\f{#1}{#2}\Big)}

\def\xp{\langle x\rangle_p}
\def\qp#1{q_p(#1)}
\def\s3{(-1)^{[\frac p3]}}
\def\s2{(-1)^{\frac{p-1}2}}

\begin{document}
 \centerline {\bf
Weighted bilinear identities and supercongruences for Ap\'ery-like polynomials}
\par\q\newline
\centerline{Yu-Tian Li$^1$ and Zhi-Hong Sun$^2$}
\par\q
\newline\centerline {$^1$School of Mathematics and Statistics} \centerline {Nanfang College}
\centerline{Guangzhou, 510970, Guangdong, China}
\centerline{Email: liyt@nfu.edu.cn}
\par\q
\newline \centerline{$^2$School of Mathematics
and Statistics} \centerline{Huaiyin Normal University}
\centerline{Huaian, Jiangsu 223300, P.R. China} \centerline{Email:
zhsun@hytc.edu.cn} \centerline{URL:
https://maths.hytc.edu.cn/szh.htm}
\vskip0.2cm
\par\q
\par {\bf Abstract.}
We establish weighted bilinear summation identities for two families of
Ap\'ery-like polynomials $g_n(x)$ and $v_n(x)$. The identities express weighted
sums in terms of consecutive endpoint values and, when necessary, lower moments.
For $g_n(x)^2$ we obtain identities with weights $(2n+1)^r$ for $1\le r\le4$;
for $v_n(x)^2$ we treat the cubic and quintic weights. Combining these formulas
with congruences for the endpoint values gives supercongruences modulo $p^3$
and $p^4$, together with special evaluations modulo $p^5$ and $p^7$, where $p$ is a prime greater than $3$. In
particular,
$$\sum_{n=0}^{p-1}(2n+1)^3v_n\!\left(\frac52\right)^2
\equiv 6p^4-\frac{143}{3}p^6\pmod {p^7},$$
confirming a congruence conjectured by Sun. The proofs use explicit quadratic
telescoping identities and $p$-adic endpoint expansions.
\par
\par\q
\newline MSC(2020): Primary 11A07, Secondary 05A19, 11B68, 11B83, 33C45
 \newline Keywords: congruence; identity; Ap\'ery-like polynomial; Bernoulli polynomial; Euler number

\section*{1. Introduction}
\par In [10], the second author investigated the properties of Ap\'ery-like polynomials $g_n(x)$ and $v_n(x)$ defined by
\begin{align*}
&g_{-1}(x)=0,\ g_0(x)=1,\ (n+1)^2g_{n+1}(x)=\Big(2n(n+1)+\f{x+1}2\Big)g_n(x)-n^2g_{n-1}(x)\ (n\ge 0),
\\&v_{-1}(x)=0,\ v_0(x)=1,\ (n+1)^3v_{n+1}(x)=(2n+1)(n(n+1)+x)v_n(x)-n^3v_{n-1}(x)\ (n\ge 0).\end{align*}
Let $G_n(x)$ and $V_n(x)$ be given by
\begin{align*}
&G_n(x)=\sum_{k=0}^n\b nk(-1)^k\b xk\b{-1-x}k,
\\&V_n(x)=\sum_{k=0}^n\b nk\b{n+k}k(-1)^k\b
xk\b{-1-x}k\ (n=0,1,2,\ldots),\end{align*}
 where $\b xk$ is the generalized
binomial coefficient given by
$$\b x0=1\qtq{and}\b xk=\f{x(x-1)\cdots(x-k+1)}{k!}\q(k\ge 1).$$
Then $$G_n(x)=g_n(2x(x+1)+1)\q\t{and}\q  V_n(x)=v_n(2x(x+1)+1)\eqno{(1.1)}$$ by
 [7,8,9].
 \par Let $p$ be an odd prime.
 In 2016, Long, Osburn and Swisher [5] proved the conjecture of Kimoto and Wakayama:
 $$\sum_{n=0}^{p-1}G_n\Big(-\f 12\Big)^2\e (-1)^{\f{p-1}2}\mod {p^3}.$$ In 2018, Liu [3] proved  the conjecture of Z.W. Sun:
 $$\sum_{n=0}^{p-1}G_n\Big(-\f 12\Big)^2\e (-1)^{\f{p-1}2}(1-7p^3B_{p-3})\mod {p^4},$$
 where $\{B_n\}$ are the Bernoulli numbers given by $B_0=1$ and $\sum_{k=0}^{n-1}\b nkB_k=0\ (n\ge 2)$.
 \par Let $\Bbb Z_{(p)}$ be the set of those rational numbers whose denominators are not divisible by $p$. Suppose that $p$ is a prime greater than $3$, $x\in\Bbb Z_{(p)}$ and $x\not\e -\f 12  \mod p$. Let $\xp$ be the least non-negative residue of $x$ modulo $p$. That is, $\xp\in\{0,1,\ldots,p-1\}$ and $x\e \xp\mod p$. In 2017, Z.W. Sun [11] conjectured that
$$\sum_{n=0}^{p-1}G_n(x)^2\e (-1)^{\xp}\f{1+2(x-\xp)/p}{1+2x}p\mod
{p^3},\eqno{(1.2)}$$
and proved the congruence modulo $p^2$. (1.2) was finally proved by Wang and Wang [12]. In [10], the second author obtained the general congruence for $\sum_{n=0}^{p-1}g_n(x)^2$ modulo $p^4$ under the assumption $\ls{2x-1}p=1$, where $\ls ap$ is the Legendre symbol.
 In [11], Z.W. Sun also proved the congruence
 $$\sum_{n=0}^{p-1}(2n+1)G_n(x)^2\e 0\mod {p^2}.$$
 Then Guo [1] showed that
 $$\sum_{n=0}^{p-1}(2n+1)G_n(x)^2\e p^2\sum_{k=0}^{p-1}\sum_{j=0}^k
 \f{(-1)^k}{k+1}\b{x+k}{2k}\b xj\b{x+j}j\b{2k}{j+k}\mod {p^4}, $$ and Liu [2] used Guo's work to confirm Sun's conjectures for $\sum_{n=0}^{p-1}(2n+1)G_n(-\f 1m)^2$ modulo $p^4$, where $m\in\{2,3,4,6\}$.
 \par
 In [10], the second author obtained the closed formulas for $$\sum_{n=0}^{p-1}g_n(x)g_n(y),\q \sum_{n=0}^{p-1}g_n(x)^2,
\ \sum_{n=0}^{p-1}(2n+1)v_n(x)v_n(y)\ \t{and}\ \sum_{n=0}^{p-1}(2n+1)v_n(x)^2$$ and gave applications to supercongruences. In particular, he showed that for $x\in\Bbb Z_{(p)}$ with $\ls{2x-1}p=1$,
$$\sum_{n=0}^{p-1}(2n+1)v_n(x)^2\e
\f{x_1+4x_0-2}{2(2x-1)}p^2+\f{x_1(x_1+8x_0-4)}{8(2x-1)^2}p^3
\mod {p^4},$$
where $x_0$ and $x_1$ are given by $x_0\in\{1,2,\ldots,\f{p-1}2\}$, $2x-1\e (2x_0-1)^2\mod p$ and $2x-1=(2x_0-1)^2+px_1$. This yields
$$\sum_{n=0}^{p-1}(2n+1)V_n(x)^2\e \f{1+2x'}{1+2x}p^2\mod
{p^4}\q\t{for $x\not\e -\f 12\mod p$}, \eqno{(1.3)}$$
 where $x'=(x-\xp)/p$.

\par In Section 2 of this paper, we establish the formulas for
\begin{align*}
&\sum_{n=0}^{p-1}(2n+1)g_n(x)g_n(y),\q \sum_{n=0}^{p-1}(2n+1)^2g_n(x)g_n(y),
\ \sum_{n=0}^{p-1}(2n+1)^3v_n(x)v_n(y),
\\&  \sum_{n=0}^{p-1}(2n+1)^3v_n(x)^2,\  \sum_{n=0}^{p-1}(2n+1)^5v_n(x)v_n(y),\ \sum_{n=0}^{p-1}(2n+1)^5v_n(x)^2,
\\&\sum_{n=0}^{p-1}(2n+1)^3g_n(x)^2\q\t{and}\q \sum_{n=0}^{p-1}(2n+1)^4g_n(x)^2.\end{align*}
In Section 3, we use the identities in Section 2 to deduce new congruences modulo prime powers. Suppose that $p>3$ is a prime and $x\in\Bbb Z_{(p)}$. We establish the congruences for
\begin{align*}&\sum_{n=0}^{p-1}(2n+1)^3v_n(x)^2,\ \sum_{n=0}^{p-1}(2n+1)^5v_n(x)^2 \mod {p^4},
\\&\sum_{n=0}^{p-1}(2n+1)V_n\Ls 12^2, \q
\sum_{n=0}^{p-1}(2n+1)^3V_n\Ls 12^2\mod {p^7},
\\&\sum_{n=0}^{p-1}(2n+1)g_n(x)^2,\q \sum_{n=0}^{p-1}(2n+1)^3g_n(x)^2\mod {p^4},\q
\\&\sum_{n=0}^{p-1}(2n+1)^2g_n(x)^2,
\q\sum_{n=0}^{p-1}(2n+1)^4g_n(x)^2\mod {p^3},
\\&\sum_{n=0}^{p-1}G_n\Ls 12^2,\q \sum_{n=0}^{p-1}G_n\Ls 32^2\mod {p^5},
\\& \sum_{n=0}^{p-1}(2n+1)G_n\Big(-\f 1m\Big)^2\ {and}\  \sum_{n=0}^{p-1}(2n+1)^3G_n\Big(-\f 1m\Big)^2\mod {p^5}\ (m=-2,2,3,4,6).
\end{align*}
As consequences, setting $x'=(x-\xp)/p$ we have
$$\sum_{n=0}^{p-1}(2n+1)G_n(x)^2\e (-1)^{\xp}(1+x'(x'+1))p^2\mod {p^4}\q\t{for $x\not\e 0,-1,-\f 12\mod p$},\eqno{(1.4)}$$
and
$$\sum_{n=0}^{p-1}(2n+1)^2G_n(x)^2
\e (-1)^{\xp-1}\f{(2x^2(x+1)^2-1)(1+2x')}{(4x^2-1)(2x+3)}p
\mod {p^3}\eqno{(1.5)}$$
for $x\not\e -\f 32,-\f 12,\f 12\mod p$. Clearly, (1.4) extends previous work of Sun, Guo and Liu.
\par In addition to the above notation, throughout this paper
we use the following notations: $\Bbb Z^+\f{\q}{\q}$the set of positive intgers, $\Bbb R\f{\q}{\q}$the set of real numbers, $[x]\f{\q}{\q}$the greatest integer not exceeding $x$, $q_p(a)=(a^{p-1}-1)/p$, $H_0=0,\ H_n=1+\f 12+\cdots+\f 1n\ (n\ge 1)$. The Bernoulli polynomials $\{B_n(x)\}$, Euler numbers $\{E_n\}$, Euler polynomials $\{E_n(x)\}$ and the sequences $\{U_n\}$ and $\{s_n\}$ are defined by
\begin{align*} &B_n(x)=\sum_{k=0}^n\b nkB_kx^{n-k}\ (n\ge 0),
\\& E_{2n-1}=0,\q E_0=1,\q E_{2n}=-\sum_{k=0}^{n-1}
\b {2n}{2k}E_{2k}\q(n\ge 1),
\\&E_n(x)=\f 1{2^n}\sum_{k=0}^n\b nk(2x-1)^{n-k}E_k\ (n\ge 0),
\\& U_{2n-1}=0,\q U_0=1,\q U_{2n}=-2\sum_{k=0}^{n-1}
\b {2n}{2k}U_{2k}\q(n\ge 1),
\\&s_0=1\q \t{and}\q s_n=1-\sum_{k=0}^{n-1} \b
nk2^{2n-1-2k}s_k\q (n\ge 1).
\end{align*}
\section*{2. New identities involving $v_n(x)$ and $g_n(x)$  }
\par\q
\par{\bf Theorem 2.1} {\sl Let $p\in\Bbb Z^+$, $x,y\in\Bbb R$, $x\not=y$ and
\[f_3(n;x,y)=-8n^3+4n^2(x-y)-4n(x+y-2)
-7x-y+4+\f{8x(x-1)}{x-y}.
\]
Then}
\begin{align*}
 &\big((x-y)^2+4(x-y)-8(x-1)\big)\sum_{n=0}^{p-1}(2n+1)^3v_n(x)v_n(y)
 \\&=8p^6\bigl(v_p(x)v_p(y)+v_{p-1}(x)v_{p-1}(y)
 \bigr)
 \\&\q+p^3f_3(p;x,y)v_p(x)v_{p-1}(y)-p^3
 f_3(-p;x,y)v_{p-1}(x)v_p(y).
\end{align*}
\par{\it Proof.} Set
\begin{align*}h_n(x,y)&= 8n^6v_n(x)v_n(y)+8n^6v_{n-1}(x)v_{n-1}(y)
  \\&\q+n^3f_3(n;x,y)v_n(x)v_{n-1}(y)-n^3f_3(-n;x,y)
 v_{n-1}(x)v_n(y).
 \end{align*}
Since $(n+1)^3v_{n+1}(x)=(2n+1)(n(n+1)+x)v_n(x)-n^3v_{n-1}(x)$, we see that
 \begin{align*}&h_{n+1}(x,y)
 \\&=8(n+1)^6v_{n+1}(x)v_{n+1}(y)+8(n+1)^6v_{n}(x)v_{n}(y)
  \\&\q+(n+1)^3f_3(n+1;x,y)v_{n+1}(x)v_{n}(y)-(n+1)^3f_3(-n-1;x,y)
 v_{n}(x)v_{n+1}(y)
 \\&=8((2n+1)(n(n+1)+x)v_n(x)-n^3v_{n-1}(x))
 ((2n+1)(n(n+1)+y)v_n(y)-n^3v_{n-1}(y))
 \\&\q+8(n+1)^6v_{n}(x)v_{n}(y)+f_3(n+1;x,y)((2n+1)(n(n+1)+x)v_n(x)-n^3v_{n-1}(x))v_n(y)
 \\&\q-f_3(-n-1;x,y) v_{n}(x)((2n+1)(n(n+1)+y)v_n(y)-n^3v_{n-1}(y))
 \\&=h_n(x,y)+((x-y)^2+4(x-y)-8(x-1))(2n+1)^3v_n(x)v_n(y).
 \end{align*}
 Thus,
 \begin{align*}&((x-y)^2+4(x-y)-8(x-1))\sum_{n=0}^{p-1}(2n+1)^3v_n(x)v_n(y)
 \\&=\sum_{n=0}^{p-1}(h_{n+1}(x,y)-h_n(x,y))=h_p(x,y)-h_0(x,y)=h_p(x,y).
 \end{align*}
 This proves the theorem.
 \vskip0.2cm

 \par{\bf Theorem 2.2} {\sl Let $p\in\Bbb Z^+$ and $x\not=1$. Then}
 \begin{align*}&\sum_{n=0}^{p-1}(2n+1)^3v_n(x)^2
 +x\sum_{n=0}^{p-1}(2n+1)v_n(x)^2
 \\&=-\f{p^6}{x-1}(v_p(x)-v_{p-1}(x))^2+2p^4v_p(x)v_{p-1}(x).
 \end{align*}
 \par{\it Proof.} From Theorem 2.1 we have
 \begin{align*}&\sum_{n=0}^{p-1}(2n+1)^3v_n(x)^2
 \\&=\lim_{y\to x}\sum_{n=0}^{p-1}(2n+1)^3v_n(x)v_n(y)
 \\&=\f{p^3}{-8(x-1)}\Big(8p^3(v_p(x)^2+v_{p-1}(x)^2)
 +(-8p^3-4p(2x-2)-8x+4)v_p(x)v_{p-1}(x)
 \\&\q-(8p^3+4p(2x-2)-8x+4)v_{p-1}(x)v_p(x)
 \\&\q+8x(x-1)\lim_{y\to x}\f{v_p(x)v_{p-1}(y)-v_{p-1}(x)v_p(y)}{x-y}\Big)
 \\&=\f{p^3}{-8(x-1)}\Big(8p^3(v_p(x)-v_{p-1}(x))^2
 -16p(x-1)v_p(x)v_{p-1}(x)
 \\&\q+8x(x-1)\lim_{y\to x}\f{-v_p(x)(v_{p-1}(x)-v_{p-1}(y))+v_{p-1}(x)(v_p(x)-v_p(y))}{x-y}\Big)
 \\&=-\f{p^6}{x-1}(v_p(x)-v_{p-1}(x))^2+2p^4v_p(x)v_{p-1}(x)
 -p^3x(v_{p-1}(x)v'_p(x)-v_p(x)v'_{p-1}(x)).
 \end{align*}
 By [10, (3.5)],
 $$\sum_{n=0}^{p-1}(2n+1)v_n(x)^2
 =p^3(v_{p-1}(x)v'_p(x)-v_p(x)v'_{p-1}(x)).$$
 Thus the result follows.
 \vskip0.2cm
 \par{\bf Theorem 2.3} {\sl Let $p\in\Bbb Z^+$, $x,y\in\Bbb R$, $x\not=y$,
\[ D_1=(x-y)^2+4(x-y)-8x+8,\quad  D_2=(x-y)^2+16(x-y)-32x+80,\]
and
\begin{align*}
 &f_5(n;x,y)\\&=-64n^5D_1+(2n+1)^4(x-y)^3
 +16n^4(4(x-y)^2-8(x-1)(x-y))
 \\ &\quad -16n^3\bigl((4x-15)(x-y)^2+(8x-32)(x-y)-8x^2+64x-104
 \bigr)\\
 &\quad+24n^2\bigl(-8(x-1)(x-y)^2+(8x^2-16x-8)(x-y)\bigr)\\
 &\quad -8n\bigl((26x-10)(x-y)^2+(-72x^2+112x+40)(x-y)
 +48x^3-112 x^2-80x+144\bigr)\\
 &\quad-20(4x-1)(x-y)^2+(464x^2-552x-56)(x-y)
 \\&\q  -768x^3
 +1536x^2+256x
 -576+\f{128(x-1)(3x^3-5x^2-7x+2)}{x-y} .
\end{align*}
Then for $D_1D_2\not=0$,}
\begin{align*}
 &\sum_{n=0}^{p-1}(2n+1)^5v_n(x)v_n(y)
 \\&=
 \frac{p^3}{D_1D_2}\Big(16p^3(4p(p-1)D_1+5(x-y)^2-(24x-32)(x-y)+24x^2-64x-8)
 v_p(x)v_p(y)
 \\&\q+16p^3(4p(p+1)D_1+5(x-y)^2-(24x-32)(x-y)+24x^2-64x-8)v_{p-1}(x)
 v_{p-1}(y)
 \\ &\quad+f_5(p;x,y)\,v_p(x)v_{p-1}(y)
 -f_5(-p;x,y)\,v_{p-1}(x)v_p(y)\Big).
\end{align*}
\par{\it Proof.} Set
$$E=5(x-y)^2-(24x-32)(x-y)+24x^2-64x-8$$
and
\begin{align*}h_n(x,y)&=
\frac{n^3}{D_1D_2}\Big(16n^3(4n(n-1)D_1+E)
 v_n(x)v_n(y)
 \\&\q+16n^3(4n(n+1)D_1+E)v_{n-1}(x)
 v_{n-1}(y)
 \\ &\quad+f_5(n;x,y)\,v_n(x)v_{n-1}(y)
 -f_5(-n;x,y)\,v_{n-1}(x)v_n(y)\Big).
\end{align*}
Using the fact that $(n+1)^3v_{n+1}(x)=(2n+1)(n(n+1)+x)v_n(x)-n^3v_{n-1}(x)$
we deduce that
\begin{align*}
h_{n+1}(x,y)&
=\frac{1}{D_1D_2}\Big(16(4n(n+1)D_1
+E) (n+1)^3v_{n+1}(x)\cdot (n+1)^3v_{n+1}(y)
 \\&\q+16(n+1)^6(4(n+1)(n+2)D_1
 +E)v_{n}(x)
 v_{n}(y)
 \\ &\quad+f_5(n+1;x,y)(n+1)^3v_{n+1}(x)v_{n}(y)
 -f_5(-n-1;x,y)\,v_{n}(x)(n+1)^3v_{n+1}(y)\Big)
  \\&
=\frac{1}{D_1D_2}\Big(16(4n(n+1)D_1
+E) ((2n+1)(n(n+1)+x)v_n(x)-n^3v_{n-1}(x))\cdot \\&\q\times((2n+1)(n(n+1)+y)v_n(y)-n^3v_{n-1}(y))
 \\&\q+16(n+1)^6(4(n+1)(n+2)D_1
 +E)v_{n}(x)
 v_{n}(y)
 \\&\quad+f_5(n+1;x,y)((2n+1)(n(n+1)+x)v_n(x)-n^3v_{n-1}(x))v_{n}(y)
 \\&\q-f_5(-n-1;x,y)\,v_{n}(x)((2n+1)(n(n+1)+y)v_n(y)-n^3v_{n-1}(y))\Big)
 \\&=h_n(x,y)+(2n+1)^5v_n(x)v_n(y).
\end{align*}
Thus,
$$\sum_{n=0}^{p-1}(2n+1)^5v_n(x)v_n(y)
=\sum_{n=0}^{p-1}(h_{n+1}(x,y)-h_n(x,y))=h_p(x,y)-h_0(x,y)=h_p(x,y).$$
The proof is now complete.
\vskip0.2cm
\par{\bf Theorem 2.4} {\sl For $p\in\Bbb Z^+$,}
\begin{align*}
& (1-x)(2x-5)\sum_{n=0}^{p-1}(2n+1)^5v_n(x)^2
-(1-x)(3x^3-5x^2-7x+2)\sum_{n=0}^{p-1}(2n+1)v_n(x)^2
\\&=p^6(4p(p-1)(x-1)-3x^2+8x+1)v_p(x)^2
+p^6(4p(p+1)(x-1)-3x^2+8x+1)v_{p-1}(x)^2
\\&\q+2p^4\bigl(4p^4(1-x)+p^2(-x^2+8x-13)+3x^3-7x^2-5x+9
 \bigr)v_p(x)v_{p-1}(x).
 \end{align*}
\par{\it Proof.} Set
\begin{align*}
h_n(x)&=n^6(4n(n-1)(x-1)-3x^2+8x+1)v_n(x)^2
\\&\q+n^6(4n(n+1)(x-1)-3x^2+8x+1)v_{n-1}(x)^2
\\&\q+2n^4\big(4n^4(1-x)+n^2(-x^2+8x-13)+3x^3-7x^2-5x+9
 \big)v_n(x)v_{n-1}(x).
 \end{align*}
 Since $(n+1)^3v_{n+1}(x)=(2n+1)(n(n+1)+x)v_n(x)-n^3v_{n-1}(x)$,
  we see that
 \begin{align*}&h_{n+1}(x)\\&=(n+1)^6(4n(n+1)(x-1)-3x^2+8x+1)v_{n+1}(x)^2
 \\&\q+(n+1)^6(4(n+1)(n+2)(x-1)-3x^2+8x+1)v_{n}(x)^2
\\&\q+2(n+1)^4\big(4(n+1)^4(1-x)+(n+1)^2(-x^2+8x-13)
\\&\q+3x^3-7x^2-5x+9
 \big)v_{n+1}(x)v_{n}(x)
 \\&=(4n(n+1)(x-1)-3x^2+8x+1)((2n+1)(n(n+1)+x)v_n(x)-n^3v_{n-1}(x))^2
 \\&\q+(n+1)^6(4(n+1)(n+2)(x-1)-3x^2+8x+1)v_{n}(x)^2
\\&\q+2(n+1)\big(4(n+1)^4(1-x)+(n+1)^2(-x^2+8x-13)
\\&\q+3x^3-7x^2-5x+9
 \big)((2n+1)(n(n+1)+x)v_n(x)-n^3v_{n-1}(x))v_n(x)
 \\&=h_n(x)+ ((1-x)(2x-5)(2n+1)^5-(1-x)(3x^3-5x^2-7x+2)
 (2n+1))v_n(x)^2.
 \end{align*}
 Thus,
 \begin{align*}
 & (1-x)(2x-5)\sum_{n=0}^{p-1}(2n+1)^5v_n(x)^2
-(1-x)(3x^3-5x^2-7x+2)\sum_{n=0}^{p-1}(2n+1)v_n(x)^2
 \\&=\sum_{n=0}^{p-1}(h_{n+1}(x)-h_n(x))=h_p(x)-h_0(x)
 =h_p(x).
 \end{align*}
 This proves the theorem.
\vskip0.2cm
\par{\bf Theorem 2.5} For every positive integer $p$,
\begin{align*}&\big((x-y)^2-4(x+y)+8\big)\sum_{n=0}^{p-1}(2n+1)g_n(x)g_n(y)
\\&=
16p^4(g_p(x)g_p(y)+g_{p-1}(x)g_{p-1}(y))
\\&\q-2p^2(8p^2-2p(x-y)+x+y-2)g_p(x)g_{p-1}(y)
\\&\q-2p^2(8p^2+2p(x-y)+x+y-2)g_{p-1}(x)g_p(y).
\end{align*}
\par{\it Proof.} Set
\begin{align*}h_n(x,y)&=16n^4(g_n(x)g_n(y)+g_{n-1}(x)g_{n-1}(y))
\\&\q-2n^2(8n^2-2n(x-y)+x+y-2)g_n(x)g_{n-1}(y)
\\&\q-2n^2(8n^2+2n(x-y)+x+y-2)g_{n-1}(x)g_n(y).
\end{align*}
Then
\begin{align*}&h_{n+1}(x,y)
\\&=16(n+1)^4(g_{n+1}(x)g_{n+1}(y)+g_{n}(x)g_{n}(y))
\\&\q-2(n+1)^2(8(n+1)^2-2(n+1)(x-y)+x+y-2)g_{n+1}(x)g_{n}(y)
\\&\q-2(n+1)^2(8(n+1)^2+2(n+1)(x-y)+x+y-2)g_{n}(x)g_{n+1}(y)
\\&=16\Big(\big(2n(n+1)+\f{x+1}2\big)g_n(x)-n^2g_{n-1}(x)\Big)
\Big(\big(2n(n+1)+\f{y+1}2\big)g_n(y)-n^2g_{n-1}(y)\Big)
\\&\q+16(n+1)^4g_n(x)g_n(y)-2(8(n+1)^2-2(n+1)(x-y)+x+y-2)
\\&\q\times\Big(\big(2n(n+1)+\f{x+1}2\big)g_n(x)-n^2g_{n-1}(x)\Big)g_n(y)
\\&\q-2(8(n+1)^2+2(n+1)(x-y)+x+y-2)
 g_n(x)\Big(\big(2n(n+1)+\f{y+1}2\big)g_n(y)-n^2g_{n-1}(y)
\Big)
\\&=h_n(x,y)+\big((x-y)^2-4(x+y)+8\big)(2n+1)g_n(x)g_n(y).
\end{align*}
Thus,
\begin{align*}&\big((x-y)^2-4(x+y)+8\big)\sum_{n=0}^{p-1}(2n+1)g_n(x)g_n(y)
\\&=\sum_{n=0}^{p-1}(h_{n+1}(x,y)-h_n(x,y))
=h_p(x,y)-h_0(x,y)=h_p(x,y).
\end{align*}
This proves the theorem.
\vskip0.2cm
\par{\bf Corollary 2.1} For every positive integer $p$ and $x\not=1$,
\begin{align*}
&\sum_{n=0}^{p-1}(2n+1)g_n(x)^2
= p^2g_p(x)g_{p-1}(x)-\frac{2p^4}{x-1}\bigl(g_p(x)
-g_{p-1}(x)\bigr)^2.
\end{align*}
\par{\it Proof.} Taking $y=x$ in Theorem 2.5 gives the result. \vskip0.2cm
\par{\bf Theorem 2.6} {\sl
Let $p\in\Bbb Z^+$. Then }
\begin{align*}&\sum_{n=0}^{p-1}\big(\big((x-y)^2-16(x+y)+80\big)
(2n+1)^2+(x-y)^2-2(x+y-2)^2+16\big)g_n(x)g_n(y)
\\&=32p^4(2p-1)g_p(x)g_p(y)+32p^4(2p+1)g_{p-1}(x)g_{p-1}(y)
\\&\quad-8p^3\big(8p^2-(x-y)p+x+y-6\big)g_p(x)g_{p-1}(y)
\\&\quad-8p^3\big(8p^2+(x-y)p+x+y-6\big)g_{p-1}(x)g_p(y).
\end{align*}
\par{\it Proof.} Let \begin{align*}h_n(x,y)&=32n^4(2n-1)g_n(x)g_n(y)+32n^4(2n+1)
g_{n-1}(x)g_{n-1}(y)
\\&\quad-8n^3\big(8n^2-(x-y)n+x+y-6\big)g_n(x)g_{n-1}(y)
\\&\quad-8n^3\big(8n^2+(x-y)n+x+y-6\big)g_{n-1}(x)g_n(y).
\end{align*}
Then
\begin{align*}&h_{n+1}(x,y)\\&=
32(n+1)^4(2n+1)g_{n+1}(x)g_{n+1}(y)+32(n+1)^4(2n+3)
g_{n}(x)g_{n}(y)
\\&\quad-8(n+1)^3\big(8(n+1)^2-(x-y)(n+1)+x+y-6\big)
g_{n+1}(x)g_{n}(y)
\\&\quad-8(n+1)^3\big(8(n+1)^2+(x-y)(n+1)
+x+y-6\big)g_{n}(x)g_{n+1}(y)
\\&=32(2n+1)\Big(\big(2n(n+1)+\f{x+1}2\big)g_n(x)-n^2g_{n-1}(x)
\Big)\\&\q\times
\Big(\big(2n(n+1)+\f{y+1}2\big)g_n(y)-n^2g_{n-1}(y)\Big)
+32(n+1)^4(2n+3)g_{n}(x)g_{n}(y)
\\&\quad-8(n+1)\big(8(n+1)^2-(x-y)(n+1)+x+y-6\big)
\\&\q\times\Big(\big(2n(n+1)+\f{x+1}2\big)g_n(x)-n^2g_{n-1}(x)\Big)g_n(y)
\\&\quad-8(n+1)\big(8(n+1)^2+(x-y)(n+1)+x+y-6\big)g_n(x)
\\&\q\times\Big(\big(2n(n+1)+\f{y+1}2\big)g_n(y)-n^2g_{n-1}(y)
\Big)
\\&=h_n(x,y)+\big(\big((x-y)^2-16(x+y)+80\big)
(2n+1)^2+(x-y)^2-2(x+y-2)^2+16\big)g_n(x)g_n(y).
\end{align*}
Hence,
\begin{align*}&\sum_{n=0}^{p-1}\big(\big((x-y)^2-16(x+y)+80\big)
(2n+1)^2+(x-y)^2-2(x+y-2)^2+16\big)g_n(x)g_n(y)
\\&=\sum_{n=0}^{p-1}(h_{n+1}(x,y)-h_n(x,y))
=h_p(x,y)-h_0(x,y)=h_p(x,y).
\end{align*}
This proves the theorem.
\vskip0.2cm

\par{\bf Corollary 2.2} {\sl For every positive integer $p$},
\begin{align*}
&(4x-10)\sum_{n=0}^{p-1}(2n+1)^2g_n(x)^2+(x^2-2x-1)
\sum_{n=0}^{p-1}g_n(x)^2
\\&=4p^4(1-2p)g_p(x)^2-4p^4(1+2p)g_{p-1}(x)^2+4p^3
(4p^2+x-3)g_p(x)g_{p-1}(x).
\end{align*}
\par{\it Proof.} Putting $y=x$ in Theorem 2.6 gives the result.
\vskip0.2cm
\par{\bf Theorem 2.7} {\sl Let $p\in\Bbb Z^+$ and $x\in\Bbb R$. Then
$$3(x-1)(x-5)\sum_{n=0}^{p-1}(2n+1)^3g_n(x)^2=h_p(x),$$
where}
\begin{align*} h_p(x)&=2p^4\big(-4p(p-1)(x-1)+x^2-5x-4\big)g_p(x)^2
\\&\quad+2p^4\big(-4p(p+1)(x-1)+x^2-5x-4\big)g_{p-1}(x)^2
\\&\quad+p^2\big(16p^4(x-1)-4p^2(3x-11)-(x-1)(x^2-3x-10)\big)
g_p(x)g_{p-1}(x).
\end{align*}
\par{\it Proof.} Let $h_0(x)=0$. It is clear that
\begin{align*} h_{n+1}(x)&=2(n+1)^4\big(-4n(n+1)(x-1)+x^2-5x-4\big)g_{n+1}(x)^2
\\&\quad+2(n+1)^4\big(-4(n+1)(n+2)(x-1)+x^2-5x-4\big)g_{n}(x)^2
\\&\quad+(n+1)^2\big(16(n+1)^4(x-1)-4(n+1)^2(3x-11)
\\&\q-(x-1)(x^2-3x-10)\big)
g_{n+1}(x)g_{n}(x)
\\&=2\big(-4n(n+1)(x-1)+x^2-5x-4\big)\Big(\big(2n(n+1)+\f{x+1}2\big)g_n(x)
-n^2g_{n-1}(x)\Big)^2
\\&\quad+2(n+1)^4\big(-4(n+1)(n+2)(x-1)+x^2-5x-4\big)g_{n}(x)^2
\\&\quad+\big(16(n+1)^4(x-1)-4(n+1)^2(3x-11)-(x-1)(x^2-3x-10)\big)
\\&\q\times\Big(\big(2n(n+1)+\f{x+1}2\big)g_n(x)
-n^2g_{n-1}(x)\Big)g_n(x)
\\&=h_n(x)+3(x-1)(x-5)(2n+1)^3g_n(x)^2.
\end{align*}
Thus,
\begin{align*}&3(x-1)(x-5)\sum_{n=0}^{p-1}(2n+1)^3g_n(x)^2
\\&=\sum_{n=0}^{p-1}(h_{n+1}(x)-h_n(x))=h_{p}(x)-h_0(x)=h_p(x).
\end{align*}
This proves the theorem.
\vskip0.2cm
\par{\bf Corollary 2.3} {\sl For $p\in\Bbb Z^+$,}
\begin{align*}
&6x(x^2-1)(x+2)\sum_{n=0}^{p-1}(2n+1)^3G_n(x)^2
\\&=p^4(-8p(p-1)x(x+1)+4x^4+8x^3-2x^2-6x-8)G_p(x)^2
\\&\q+p^4(-8p(p+1)x(x+1)+4x^4+8x^3-2x^2-6x-8)G_{p-1}(x)^2
\\&\q+p^2(16p^4x(x+1)-4p^2(3x^2+3x-4)-2x(x^2-1)(x+2)(2x^2+2x+3))
G_p(x)G_{p-1}(x).
\end{align*}
\par{\it Proof.} Since $G_n(x)=g_n(2x^2+2x+1)$, replacing $x$ with $2x^2+2x+1$ in Theorem 2.7 yields the result.
\vskip0.2cm
\par{\bf Theorem 2.8} {\sl Let $p\in\Bbb Z^+$ and $x\in\Bbb R$. Then}
$$8(2x-5)(2x-17)\sum_{n=0}^{p-1}(2n+1)^4g_n(x)^2
-(3x^4-20x^3-78x^2+268x-5)\sum_{n=0}^{p-1}g_n(x)^2=h_p(x),$$
where
\begin{align*}
h_p(x)&=4p^4(2p-1)\big(-8p(p-1)(2x-5)+3x^2-26x-25\big)g_p(x)^2
\\&\quad+4p^4(2p+1)\big(-8p(p+1)(2x-5)+3x^2-26x-25\big)g_{p-1}(x)^2
\\&\quad+4p^3\big(32p^4(2x-5)+4p^2(x^2-28x+135)-(x+5)(3x^2-34x+71)\big)
g_p(x)g_{p-1}(x).
\end{align*}
\par{\it Proof.} Set $h_0(x)=0$. It is clear that
\begin{align*}
&h_{n+1}(x)\\&=4(n+1)^4(2n+1)\big(-8n(n+1)(2x-5)+3x^2-26x-25\big)
g_{n+1}(x)^2
\\&\quad+4(n+1)^4(2n+3)\big(-8(n+1)(n+2)(2x-5)
+3x^2-26x-25\big)g_{n}(x)^2
\\&\quad+4(n+1)^3\big(32(n+1)^4(2x-5)
+4(n+1)^2(x^2-28x+135)
\\&\q-(x+5)(3x^2-34x+71)\big)
g_{n+1}(x)g_{n}(x)
\\&=4(2n+1)\big(-8n(n+1)(2x-5)+3x^2-26x-25\big)
\\&\q\times\Big(\big(2n(n+1)+\f{x+1}2\big)g_n(x)
-n^2g_{n-1}(x)\Big)^2
\\&\q\quad+4(n+1)^4(2n+3)\big(-8(n+1)(n+2)(2x-5)
+3x^2-26x-25\big)g_{n}(x)^2
\\&\q+4(n+1)\big(32(n+1)^4(2x-5)
+4(n+1)^2(x^2-28x+135)-(x+5)(3x^2-34x+71)\big)
\\&\q\times \Big(\big(2n(n+1)+\f{x+1}2\big)g_n(x)
-n^2g_{n-1}(x)\Big)g_{n}(x)
\\&=h_n(x)+\big(8(2x-5)(2x-17)(2n+1)^4
-(3x^4-20x^3-78x^2+268x-5)\big)g_n(x)^2.
\end{align*}
Therefore,
\begin{align*}&\sum_{n=0}^{p-1}\big(8(2x-5)(2x-17)(2n+1)^4
-(3x^4-20x^3-78x^2+268x-5)\big)g_n(x)^2
\\&=\sum_{n=0}^{p-1}(h_{n+1}(x)-h_n(x))=h_p(x)-h_0(x)=h_p(x),
\end{align*}
which concludes the proof.

\section*{3. Supercongruences involving $v_n(x)$ and $g_n(x)$}

 For an odd prime $p$ and $x\in\Bbb Z_{(p)}$ with $\ls{2x-1}p=1$ let $x_0$ and $x_1$ be given by $$x_0\in\big\{1,2,\ldots,\f{p-1}2\big\},\ 2x-1\e (2x_0-1)^2\mod p,
  \ 2x-1=(2x_0-1)^2+px_1.\eqno{(3.1)}$$

\par{\bf Lemma 3.1 ([10, Theorem 4.4])} {\sl Let $p$ be an odd prime, $x\in\Bbb Z_{(p)}$ and $\ls{2x-1}p=1$. Then
$$\sum_{n=0}^{p-1}(2n+1)v_n(x)^2\e
\f{x_1+4x_0-2}{2(2x-1)}p^2+\f{x_1(x_1+8x_0-4)}{8(2x-1)^2}p^3
\mod {p^4},$$
where $x_0$ and $x_1$ are given by (3.1). }
\vskip0.2cm
\par{\bf Theorem 3.1} {\sl Let $p$ be an odd prime, $x\in\Bbb Z_{(p)}$ and $\ls{2x-1}p=1$, and let $x_0$ and $x_1$ be given by (3.1).
\par $(\t{\rm i})$ If $x\not\e 1\mod {p^3}$, then
$$\sum_{n=0}^{p-1}(2n+1)^3v_n(x)^2\e
-\f{x(x_1+4x_0-2)}{2(2x-1)}p^2-\f{xx_1(x_1+8x_0-4)}{8(2x-1)^2}p^3
\mod {p^4}.$$
\par $(\t{\rm ii})$ If $x\not\e 1,\f 52\mod p$, then}
$$\sum_{n=0}^{p-1}(2n+1)^5v_n(x)^2\e
\f{3x^3-5x^2-7x+2}{2x-5}\Big(\f{x_1+4x_0-2}{2(2x-1)}
+\f{x_1(x_1+8x_0-4)}{8(2x-1)^2}p\Big)p^2\mod {p^4}.$$
\par{\it Proof.} This is immediate from Theorem 2.2, 2.4 and Lemma 3.1.
\vskip0.2cm
\par{\bf Lemma 3.2 ([10, Corollary 4.2])} {\sl Let $p$ be an odd prime, $x\in\Bbb Z_{(p)}$ and $x\not\e -\f 12\mod p$. Then
$$\sum_{n=0}^{p-1}(2n+1)V_n(x)^2\e \f{1+2x'}{1+2x}p^2\mod
{p^4}, $$ where $x'=(x-\xp)/p$.}
\vskip0.2cm
\par{\bf Theorem 3.2} {\sl Let $p$ be an odd prime, $x\in\Bbb Z_{(p)}$, $x\not\e -\f 12\mod p$ and $x'=(x-\xp)/p$.
\par $(\t{\rm i})$ If $x(x+1)\not\e 0\mod {p^3}$, then
$$\sum_{n=0}^{p-1}(2n+1)^3V_n(x)^2\e -\f{(2x^2+2x+1)(2x'+1)}{2x+1}p^2\mod {p^4}.$$
\par $(\t{\rm ii})$ If $x\not\e 0,-1,\f 12,-\f 32\mod p$, then
$$\sum_{n=0}^{p-1}(2n+1)^5V_n(x)^2\e
\f{(24x^6+72x^5+88x^4+56x^3-16x-7)(1+2x')}{(2x+3)(4x^2-1)}p^2
\mod {p^4}.$$}
\par{\it Proof.} Since $V_n(x)=v_n(2x^2+2x+1)$, substituting $x$ with $2x^2+2x+1$ in Theorems 2.2, 2.4 and then applying Lemma 3.2 yields the results.
\vskip0.2cm
\par{\bf Lemma 3.3} {\sl Let $p>3$ be a prime. Then
\begin{align*}&V_p\Ls 12\e \f 1{2\cdot 16^{p-1}}( 1+4p+8p^2)\mod{p^3},
 \\&V_{p-1}\Ls 12\e 1-4(1-\qp 2)p+(6\qp 2^2-16\qp 2+8)p^2\mod {p^3}\end{align*}
 and so $V_p\ls 12V_{p-1}\ls 12\e \f 12\mod {p^3}$. }
 \vskip0.2cm
\par{\it Proof.} It is well known (see [6]) that
$B_1=-\f 12$, $$B_{2n+1}=0\qtq{and}B_n(-x)=(-1)^nB_n(1+x)=(-1)^n\big(B_n(x)+nx^{n-1}\big)
\q(n\ge 1).$$ Thus,
$$(-1)^n\f{B_n(-x)-B_n}n
= \f{B_n(x)-B_n}n+x^{n-1}\qtq{for}n>1.$$
By [7, (2.16)],
$$\f{B_{p^2(p-1)}(\f 12)-B_{p^2(p-1)}}{p^2(p-1)}
\e  -2\qp 2+p\qp 2^2\mod
{p^2}.$$
Hence,
$$\f{B_{p^2(p-1)}(-\f 12)-B_{p^2(p-1)}}{p^2(p-1)}
\e 2-2\qp 2+p\qp 2^2\mod {p^2}. \eqno{(3.2)}$$
From [9, proof of Theorem 3.2], for $x\in\Bbb Z_{(p)}$ we have
\begin{align*} &V_p(x)-1+2\b
xp\b{-1-x}p \\&\e 2p^2\Big(\f{B_{p^2(p-1)}(-x)-B_{p^2(p-1)}}{p^2(p-1)}\Big)^2 +2p
\f{B_{p^2(p-1)}(-x)-B_{p^2(p-1)}}{p^2(p-1)}\mod {p^3}.
\end{align*}
Taking $x=\f 12$ in the congruence and then applying (3.2) we
deduce that
\begin{align*} &V_p\Ls 12-1+2\b{1/2}p\b{-1-1/2}p
\\&\e 2p^2(2-2q_p(2)+pq_p(2)^2)^2+2p(2-2q_p(2)+pq_p(2)^2)
\\&\e 4p(1-q_p(2))+(8-16q_p(2)+10q_p(2)^2)p^2\mod {p^3}.
\end{align*}
It is easily seen that
$$\b{1/2}p\b{-1-1/2}p=\f{1+2p}{1-2p}\b{-1/2}p^2=\f{(1+2p)(1+2p+4p^2)}{1-8p^3}
\cdot\f{\b{2p}p^2}{16^p}.$$ Since $\b{2p}p=2\b{2p-1}{p-1}\e 2\mod
{p^3}$, we derive
\begin{align*} 16^pV_p\Ls 12&\e
-2(1+2p)(1+2p+4p^2)\b{2p}p^2\\&\q+16^p\big(1+4p(1-q_p(2))+(8-16q_p(2)+10q_p(2)^2)p^2
\big)
\\&\e -8(1+2p)(1+2p+4p^2)+16(1+4pq_p(2)+6p^2q_p(2)^2)
\\&\q\times
\big(1+4p(1-q_p(2))+(8-16q_p(2)+10q_p(2)^2)p^2 \big)
\\&\e  8+32p+64p^2\mod{p^3},\end{align*}
which yields the first result. From [9, proof of Theorem 3.3],
$$ V_{p-1}(x)\e 1-2p\f{B_{p^2(p-1)}(-x)-B_{p^2(p-1)}}{p^2(p-1)}
+2p^2\Big(\f{B_{p^2(p-1)}(-x)-B_{p^2(p-1)}}{p^2(p-1)}\Big)^2\mod
{p^3}.$$
Taking $x=\f 12$ and then applying (3.2) gives
\begin{align*}V_{p-1}\Ls 12&\e 1-2p\big(2-2\qp 2+p\qp 2^2\big)+2p^2(2-2\qp 2)^2
\\&=1-4(1-\qp 2)p+(6\qp 2^2-16\qp 2+8)p^2\mod {p^3}.
\end{align*}
Thus,
\begin{align*}V_p\Ls 12V_{p-1}\Ls 12&\e \f{1}{2\cdot 16^{p-1}}
(1+4p+8p^2)\big(1-4(1-\qp 2)p+(6\qp 2^2-16\qp 2+8)p^2\big)
\\&\e \f{1+4p\qp 2+6p^2\qp 2^2}{2(1+p\qp 2)^4}\e\f 12\mod {p^3}.
\end{align*}
This completes the proof.
\vskip0.2cm
\par{\bf Theorem 3.3} {\sl Let $p$ be a prime greater than $3$. Then}
$$\sum_{n=0}^{p-1}(2n+1)V_n\Ls 12^2
=\sum_{n=0}^{p-1}(2n+1)v_n\Ls 52^2\e -2p^4+19p^6\mod {p^7}$$
and
$$\sum_{n=0}^{p-1}(2n+1)^3V_n\Ls 12^2=\sum_{n=0}^{p-1}(2n+1)^3v_n\Ls 52^2\e 6p^4-\f{143}3p^6\mod {p^7}.$$
\par{\it Proof.} Since $2\sls 12^2+2\cdot \f 12+1=\f 52$ we have $V_n\sls 12=v_n(\f 52)$. Now, taking $x=\f 52$ in Theorem 2.4 and then applying Lemma 3.3 gives
\begin{align*}&\f 3{16}\sum_{n=0}^{p-1}(2n+1)V_n\Ls 12^2
\\&=p^6\Big(6p(p-1)+\f 94\Big)V_p\Ls 12^2+p^6\Big(6p(p+1)+\f 94\Big)V_{p-1}\Ls 12^2
\\&\q+2p^4\Big(-6p^4+\f 34p^2-\f 38\Big)V_p\Ls 12V_{p-1}\Ls 12
\\&\e p^6\cdot \f 94\cdot \f 14+p^6\cdot\f 94\cdot 1+\f 32p^6\cdot \f 12
-\f 34p^4\cdot \f 1{2}
= \f{57}{16}p^6-\f 38p^4
\mod {p^7},\end{align*}
which yields the first result.
\par Since $v_n\sls 52=V_n\sls 12$, taking $x=\f 52$ in Theorem 2.2 and then applying Lemma 3.3
 yields
\begin{align*}&\sum_{n=0}^{p-1}(2n+1)^3V_n\Ls 12^2
 +\f 52\sum_{n=0}^{p-1}(2n+1)V_n\Ls 12^2
 \\&=-\f{p^6}{\f 52-1}\Big(V_p\Ls 12-V_{p-1}\Ls 12\Big)^2+2p^4V_p\Ls 12V_{p-1}\Ls 12
 \\&\e -\f 23p^6\Big(\f 12-1\Big)^2+2p^4\cdot \f 12=p^4-\f {p^6}6\mod {p^7}.
 \end{align*}
 Therefore,
 $$\sum_{n=0}^{p-1}(2n+1)^3V_n\Ls 12^2
 \e p^4-\f {p^6}6-\f 52(-2p^4+19p^6)=6p^4-\f{143}3p^6\mod {p^7}.$$
 This concludes the proof.
\vskip0.2cm
\par{\bf Remark 3.1} Let $p>3$ be a prime. In [10, Corollary 4.3, Conjecture 4.2 and Remark 4.1], the second author proved $\sum_{n=0}^{p-1}(2n+1)V_n\sls 12^2\e 0\mod {p^4}$ and conjectured Theorem 3.3.
\vskip0.2cm
\par{\bf Lemma 3.4 ([10, Lemma 5.3])} {\sl Let $p$ be an odd prime and $x\in\Bbb Z_{(p)}$ with $\ls{2x-1}p=1$. Then
\begin{align*}&g_p(x)\e
 \Big(1+\f{x_1(x_1+8x_0-4)}{16(2x-1)}\Big)(1+2pH_{x_0-1})
+p\f{x_1^2(x_1+12x_0-6)}{32(2x-1)^2}
\mod {p^2},
\\&g_{p-1}(x)\e
(-1)^{x_0-1}(1-2pH_{x_0-1})
\mod {p^2},
\end{align*}
where $x_0$ and $x_1$ are given by (3.1). }
\vskip0.2cm
\par{\bf Theorem 3.4} {\sl Let $p$ be an odd prime, $x\in\Bbb Z_{(p)}$, $x\not\e 1\mod p$ and $\sls{2x-1}p=1$. Then
\begin{align*}&\sum_{n=0}^{p-1}(2n+1)g_n(x)^2
\\&\e (-1)^{x_0-1}\Big(1+\f{x_1(x_1+8x_0-4)}{16(2x-1)}\Big)p^2
+(-1)^{x_0-1}\f{x_1^2(x_1+12x_0-6)}{32(2x-1)^2}p^3\mod {p^4},
\end{align*}
where $x_0$ and $x_1$ are given by (3.1). }
\vskip0.2cm
\par{\it Proof.} By Corollary 2.1,
$$\sum_{n=0}^{p-1}(2n+1)g_n(x)^2\e p^2g_p(x)g_{p-1}(x)\mod {p^4}.$$
From Lemma 3.4,
$$g_p(x)g_{p-1}(x)\e (-1)^{x_0-1}\Big(1+\f{x_1(x_1+8x_0-4)}{16(2x-1)}
+\f{x_1^2(x_1+12x_0-6)}{32(2x-1)^2}p\Big)\mod {p^2}.\eqno{(3.3)}$$
Thus the result follows.
\vskip0.2cm

\par{\bf Theorem 3.5} {\sl Let $p$ be an odd prime, $x\in\Bbb Z_{(p)}$, $x\not\e 0,-1,-\f 12\mod p$ and $x'=(x-\xp)/p$.
Then
$$\sum_{n=0}^{p-1}(2n+1)G_n(x)^2\e (-1)^{\xp}(1+x'(x'+1))p^2\mod {p^4}.$$}
\par{\it Proof}. Observe that
$$2(2x^2+2x+1)-1=(2x+1)^2\e (2(\xp+1)-1)^2
\e (2(p-\xp)-1)^2\mod p.$$
Set
$$ x_0=\begin{cases}\xp+1=x+1-px'&\t{if $\xp<\f{p-1}2$,}
\\p-\xp=(x'+1)p-x&\t{if $\xp\ge \f{p-1}2$}
\end{cases}$$
and
$$ x_1=\f{(2x+1)^2-(2x_0-1)^2}p=\begin{cases}
4x'(2x+1-px')&\t{if $\xp<\f{p-1}2$,}
\\ 4(x'+1)(2x+1-p(x'+1))&\t{if $\xp\ge \f{p-1}2$.}
\end{cases}$$
Then $x_0\in\{1,2,\ldots,\f{p-1}2\}$ and
$2(2x^2+2x+1)-1=(2x_0-1)^2+px_1$.
\par We first assume that $\xp<\f{p-1}2$. Then
\begin{align*}x_1(x_1+8x_0-4)
&=4x'(2x+1-px')(4x'(2x+1-px')+8(x+1-px')-4)
\\&\e 16(2x+1)x'\big((2x+1)(x'+1)-px'(2x'+3)\big)
\mod {p^2}
\end{align*}
and
\begin{align*}x_1^2(x_1+12x_0-6)&\e (4x'(2x+1))^2(4x'(2x+1)+12(x+1)-6)
\\&=32(2x+1)^3{x'}^2(2x'+3)\mod p.
\end{align*}
Since $2(2x^2+2x+1)-1=(2x+1)^2$ and $G_n(x)=g_n(2x^2+2x+1)$, from the above and Theorem 3.4 we deduce that
\begin{align*}&\sum_{n=0}^{p-1}(2n+1)G_n(x)^2
\\&\e (-1)^{x_0-1}\Big(1+\f{x_1(x_1+8x_0-4)}{16(2x+1)^2}\Big)p^2
+(-1)^{x_0-1}\f{x_1^2(x_1+12x_0-6)}{32(2x+1)^4}p^3
\\&\e (-1)^{\xp}\Big(1+\f{x'
\big((2x+1)(x'+1)-px'(2x'+3)\big)}{2x+1}\Big)p^2
+(-1)^{\xp}\f{{x'}^2(2x'+3)}{2x+1}p^3
\\&=(-1)^{\xp}(1+x'(x'+1))p^2\mod {p^4}.
\end{align*}
\par Now assume that $\xp\ge \f{p-1}2$. Then
\begin{align*}&x_1(x_1+8x_0-4)
\\&=4(x'+1)(2x+1-p(x'+1))\big(4(x'+1)(2x+1-p(x'+1))+8((x'+1)p-x)-4\big)
\\&=16(x'+1)(2x+1-p(x'+1))((2x+1)x'-p(x'+1)(x'-1))
\\&\e 16(2x+1)(x'+1)((2x+1)x'-(x'+1)(2x'-1)p)\mod {p^2}
\end{align*}
and
\begin{align*}x_1^2(x_1+12x_0-6)
&\e (4(x'+1)(2x+1))^2(4(x'+1)(2x+1)+12(-x)-6)
\\&=32(2x+1)^3(x'+1)^2(2x'-1)\mod p.
\end{align*}
Since $2(2x^2+2x+1)-1=(2x+1)^2$ and $G_n(x)=g_n(2x^2+2x+1)$, from the above and Theorem 3.4 we deduce that
\begin{align*}&\sum_{n=0}^{p-1}(2n+1)G_n(x)^2
\\&\e (-1)^{x_0-1}\Big(1+\f{x_1(x_1+8x_0-4)}{16(2x+1)^2}\Big)p^2
+(-1)^{x_0-1}\f{x_1^2(x_1+12x_0-6)}{32(2x+1)^4}p^3
\\&\e (-1)^{p-1-\xp}\Big(1+\f{(x'+1)((2x+1)x'-(x'+1)(2x'-1)p)}
{2x+1}\Big)p^2
\\&\q+(-1)^{p-1-\xp}\f{(x'+1)^2(2x'-1)}{2x+1}p^3
\\&=(-1)^{\xp}(1+x'(x'+1))p^2\mod {p^4}.
\end{align*}
This completes the proof.
\vskip0.2cm
\par{\bf Lemma 3.5 ([10, Theorem 5.6])} {\sl Let $p$ be an odd prime, $x\in\Bbb Z_{(p)}$ and $\ls{2x-1}p=1$. Then
\begin{align*}\sum_{n=0}^{p-1}g_n(x)^2&\e
(-1)^{x_0-1}\f{x_1+4x_0-2}{2(2x-1)}p
+(-1)^{x_0-1}\f{x_1(x_1+8x_0-4)}{8(2x-1)^2}p^2
\mod {p^3},\end{align*}
where $x_0$ and $x_1$ are given by (3.1). }
\vskip0.2cm
\par{\bf Theorem 3.6} {\sl Let $p$ be an odd prime, $x\in\Bbb Z_{(p)}$, $x\not\e \f 52\mod p$ and $\ls{2x-1}p=1$. Then
\begin{align*}\sum_{n=0}^{p-1}(2n+1)^2g_n(x)^2&\e
(-1)^{x_0}\f{x^2-2x-1}{4x-10}\Big(\f{x_1+4x_0-2}{2(2x-1)}p
+\f{x_1(x_1+8x_0-4)}{8(2x-1)^2}p^2\Big)
\mod {p^3},\end{align*}
where $x_0$ and $x_1$ are given by (3.1). }
\vskip0.2cm
\par{\it Proof.} By Corollary 2.2,
$$\sum_{n=0}^{p-1}(2n+1)^2g_n(x)^2\e -\f{x^2-2x-1}{4x-10}\sum_{n=0}^{p-1}g_n(x)^2\mod {p^3}.\eqno{(3.4)}$$
This together with Lemma 3.5 yields the result.
\vskip0.2cm
\par{\bf Theorem 3.7} {\sl Let $p$ be an odd prime, $x\in\Bbb Z_{(p)}$, $x\not\e -\f 32,-\f 12,\f 12\mod p$ and $x'=(x-\xp)/p$. Then
$$\sum_{n=0}^{p-1}(2n+1)^2G_n(x)^2
\e (-1)^{\xp-1}\f{(2x^2(x+1)^2-1)(1+2x')}{(4x^2-1)(2x+3)}p
\mod {p^3}.$$}
\par{\it Proof.} Since $G_n(x)=g_n(2x^2+2x+1)$
and
$$\f{(2x^2+2x+1)^2-2(2x^2+2x+1)-1}{4(2x^2+2x+1)-10}
=\f{2x^2(x+1)^2-1}{(2x-1)(2x+3)},$$ from (3.4) we have
$$\sum_{n=0}^{p-1}(2n+1)^2G_n(x)^2\e -\f{2x^2(x+1)^2-1}{(2x-1)(2x+3)}\sum_{n=0}^{p-1}G_n(x)^2\mod {p^3}.$$
Now, applying (1.2) yields the result.
\vskip0.2cm
\par{\bf Theorem 3.8} {\sl Let $p$ be an odd prime. Then}
$$\sum_{n=0}^{p-1}g_n\Ls 52^2=\sum_{n=0}^{p-1}G_n\Ls 12^2
=\sum_{n=0}^{p-1}G_n\Big(-\f 32\Big)^2
\e 6(-1)^{\f{p-1}2}p^3-7p^4\mod {p^5}.$$
\par{\it Proof.} For $x=\f 52$ we have $x_0=\f{p-1}2$ and $x_1=4-p$ by (3.1). Thus,
$$\f{x_1(x_1+8x_0-4)}{16(2x-1)}=\f{(4-p)(4-p+4(p-1)-4)}{16(5-1)}
\e \f{p-1}4\mod {p^2}$$
and
$$\f{x_1^2(x_1+12x_0-6)}{32(2x-1)^2}\e \f{4^2(4-6-6)}{32(5-1)^2}=-\f 14\mod p.$$
This together with Lemma 3.4 and (3.3) yields
$$g_p\Ls 52\e 1+\f{p-1}4\e \f 34\mod p,\q
g_{p-1}\Ls 52\e (-1)^{\f{p+1}2}\mod p$$
and
 $$g_p\Ls 52g_{p-1}\Ls 52\e (-1)^{\f{p-3}2}\Big(1+\f{p-1}4\Big)+(-1)^{\f{p-3}2}\Big(-\f p4\Big)=\f 34(-1)^{\f{p+1}2}\mod {p^2}.$$
Now, taking $x=\f 52$ in Corollary 2.2 and then applying the above gives
\begin{align*}
\sum_{n=0}^{p-1}g_n\Ls 52^2
&=16p^4(1-2p)g_p\Ls 52^2-16p^4(1+2p)g_{p-1}\Ls 52^2
+16p^3\Big(4p^2-\f 12\Big)g_p\Ls 52g_{p-1}\Ls 52
\\&\e 16p^4\cdot \f 9{16}-16p^4\cdot 1+16p^3\Big(-\f 12\Big)\cdot\f 34(-1)^{\f{p+1}2}=6(-1)^{\f{p-1}2}p^3-7p^4\mod {p^5}.\end{align*}
To complete the proof, we note that $G_n(\f 12)=G_n(-\f 32)=g_n(\f 52)$ since $G_n(x)=g_n(2x^2+2x+1)$.
\vskip0.2cm
\par{\bf Lemma 3.6 ([7, Theorem 6.2])} Let
$p>3$ be a prime. Then
\begin{align*} &16^pG_p\Big(-\f 12\Big)\e 12+64(-1)^{\f{p-1}2}p^2E_{p-3}\mod {p^3},
\\& 27^pG_p\Big(-\f 13\Big)\e 21+243(-1)^{[\f p3]}p^2U_{p-3}\mod {p^3},
\\&64^pG_p\Big(-\f 14\Big)\e 52+1024(-1)^{[\f p4]}p^2s_{p-3}\mod {p^3},
\\& 432^pG_p\Big(-\f 16\Big)\e 372+8640(-1)^{\f{p-1}2}p^2E_{p-3}\mod {p^3}.
\end{align*}
\par{\bf Lemma 3.7 ([4])} {\sl Let $p>3$ be a prime. Then}
$$16^{p-1}G_{p-1}\Big(-\f 12\Big)\e (-1)^{\f{p-1}2}256^{p-1}+3p^2E_{p-3}\mod{p^3}.$$

\par{\bf Lemma 3.8 ([7, Theorem 6.1])} {\sl Let $p>3$ be a prime. Then}
 \begin{align*}&27^{p-1}G_{p-1}\Big(-\f 13\Big)\e (-1)^{[\f p3]}729^{p-1}+7p^2U_{p-3}\mod {p^3},
\\&432^{p-1}G_{p-1}\Big(-\f 16\Big)\e (-1)^{\f{p-1}2}186624^{p-1}+\f{155}9p^2E_{p-3}\mod
{p^3},
\\&64^{p-1}G_{p-1}\Big(-\f 14\Big)\e (-1)^{[\f p4]}4096^{p-1}+13p^2s_{p-3}\mod {p^3}.\end{align*}

\par{\bf Lemma 3.9 } {\sl Let $p>3$ be a prime. Then}
 \begin{align*}&G_p\Big(-\f 12\Big)G_{p-1}\Big(-\f 12\Big)
 \e \f 34(-1)^{\f{p-1}2}+\f{25}4p^2E_{p-3}\mod {p^3},
 \\&G_p\Big(-\f 13\Big)G_{p-1}\Big(-\f 13\Big)
 \e \f 79(-1)^{[\f p3]}+\f{130}9p^2U_{p-3}\mod {p^3},
 \\&G_p\Big(-\f 14\Big)G_{p-1}\Big(-\f 14\Big)
 \e \f{13}{16}(-1)^{[\f p4]}+\f{425}{16}p^2s_{p-3}\mod {p^3},
 \\&G_p\Big(-\f 16\Big)G_{p-1}\Big(-\f 16\Big)
 \e \f {31}{36}(-1)^{\f{p-1}2}+\f{11285}{324}p^2E_{p-3}\mod {p^3}.
 \end{align*}
 \par{\it Proof.} From Lemmas 3.6 and 3.7,
 \begin{align*}G_p\Big(-\f 12\Big)G_{p-1}\Big(-\f 12\Big)
 &\e \f 1{16^p}\big(12+64(-1)^{\f{p-1}2}p^2E_{p-3}\big)
 \cdot\f 1{16^{p-1}}\big((-1)^{\f{p-1}2}256^{p-1}+3p^2E_{p-3}\big)
 \\&\e \f 1{16\cdot 256^{p-1}}\big(12(-1)^{\f{p-1}2}256^{p-1}+36p^2E_{p-3}+64p^2E_{p-3}
 \big)
 \\&\e \f 34(-1)^{\f{p-1}2}+\f{25}4p^2E_{p-3}\mod {p^3}.
 \end{align*}
 The remaining parts can be proved similarly.
 \vskip0.2cm
 \par{\bf Theorem 3.9} {\sl Let $p>5$ be a prime. Then
 \begin{align*}&\sum_{n=0}^{p-1}(2n+1)G_n\Big(-\f 12\Big)^2
 \e \f 34\s2 p^2+\Big(-6\s2+\f{25}4(1+E_{p-3})\Big)p^4\mod {p^5},
\\&\sum_{n=0}^{p-1}(2n+1)G_n\big(-\f 13\big)^2
\e\f 79(-1)^{[\f p3]}p^2+\Big(\f {65}9-7(-1)^{[\f p3]}+\f{130}9U_{p-3}\Big)p^4\mod {p^5},
\\&\sum_{n=0}^{p-1}(2n+1)G_n\big(-\f 14\big)^2
\e\f{13}{16}(-1)^{[\f p4]}p^2+\Big(\f{425}{48}-\f{26}3(-1)^{[\f p4]}+\f{425}{16}s_{p-3}\Big)p^4\mod {p^5}
\end{align*}
and
\begin{align*}
\sum_{n=0}^{p-1}(2n+1)G_n\big(-\f 16\big)^2
\e\f {31}{36}(-1)^{\f{p-1}2}p^2+\Big(\f{2257}{180}-\f{62}5
\s2 +\f{11285}{324}E_{p-3}\Big)p^4\mod {p^5}.
\end{align*}}
\par{\it Proof.} Since $G_n(x)=g_n(2x^2+2x+1)$, from Corollary 2.1 we see that for $x\not=0,-1$,
$$\sum_{n=0}^{p-1}(2n+1)G_n(x)^2=p^2G_p(x)G_{p-1}(x)-\f{p^4}{x(x+1)}
\big(G_p(x)-G_{p-1}(x)\big)^2.\eqno{(3.5)}$$
This together with Lemmas 3.6-3.9 yields
\begin{align*}
&\sum_{n=0}^{p-1}(2n+1)G_n\big(-\f 12\big)^2
\\&=p^2G_p\big(-\f 12\big)G_{p-1}\big(-\f 12\big)
-\f{p^4}{-\f 12\cdot \f 12}\Big(G_p\big(-\f 12\big)-G_{p-1}\big(-\f 12\big)\Big)^2
\\&\e p^2\Big(\f 34\s2+\f{25}4p^2E_{p-3}\Big)+4p^4\Big(\f{3}{4}
-(-1)^{\f{p-1}2}\Big)^2
\\&=\f 34\s2 p^2+\Big(-6\s2+\f{25}4(1+E_{p-3})\Big)p^4\mod {p^5},
\\&\sum_{n=0}^{p-1}(2n+1)G_n\big(-\f 13\big)^2
\\&=p^2G_p\big(-\f 13\big)G_{p-1}\big(-\f 13\big)
-\f{p^4}{-\f 13\cdot \f 23}\Big(G_p\big(-\f 13\big)-G_{p-1}\big(-\f 13\big)\Big)^2
\\&\e p^2\Big(\f 79(-1)^{[\f p3]}+\f{130}9p^2U_{p-3}\Big)
+\f 92p^4\Big(\f 79-(-1)^{[\f p3]}\Big)^2
\\&=\f 79(-1)^{[\f p3]}p^2+\Big(\f {65}9-7(-1)^{[\f p3]}+\f{130}9U_{p-3}\Big)p^4\mod {p^5},
\\&\sum_{n=0}^{p-1}(2n+1)G_n\big(-\f 14\big)^2
\\&=p^2G_p\big(-\f 14\big)G_{p-1}\big(-\f 14\big)
-\f{p^4}{-\f 14\cdot \f 34}\Big(G_p\big(-\f 14\big)-G_{p-1}\big(-\f 14\big)\Big)^2
\\&\e p^2\Big(\f{13}{16}(-1)^{[\f p4]}+\f{425}{16}p^2s_{p-3}\Big)+\f{16}3p^4\Big(\f{13}{16}-(-1)^{[\f p4]}\Big)^2
\\&=\f{13}{16}(-1)^{[\f p4]}p^2+\Big(\f{425}{48}-\f{26}3(-1)^{[\f p4]}+\f{425}{16}s_{p-3}\Big)p^4\mod {p^5}
\end{align*}
and
\begin{align*}
\sum_{n=0}^{p-1}(2n+1)G_n\big(-\f 16\big)^2
&=p^2G_p\big(-\f 16\big)G_{p-1}\big(-\f 16\big)
-\f{p^4}{-\f 16\cdot \f 56}\Big(G_p\big(-\f 16\big)-G_{p-1}\big(-\f 16\big)\Big)^2
\\&\e p^2\Big(\f {31}{36}(-1)^{\f{p-1}2}+\f{11285}{324}p^2E_{p-3}\Big)
+\f{36}5p^4\Big(\f{31}{36}-(-1)^{\f{p-1}2}\Big)^2
\\&=\f {31}{36}(-1)^{\f{p-1}2}p^2+\Big(\f{2257}{180}-\f{62}5
\s2 +\f{11285}{324}E_{p-3}\Big)p^4\mod {p^5}.
\end{align*}
This proves the theorem.
\vskip0.2cm
\par{\bf Theorem 3.10} {\sl Let $p>3$ be a prime. Then}
$$\sum_{n=0}^{p-1}(2n+1)G_n\Ls 12^2\e
\f 34(-1)^{\f{p+1}2}p^2+\Big(\f{125}{12}+2(-1)^{\f{p+1}2}-\f{25}4
E_{p-3}\Big)p^4\mod {p^5}.$$
\par{\it Proof.}
Since $E_n(x)+E_n(x+1)=2x^n$ and $E_n\ls 12=\f{E_n}{2^n}$ we see that
$E_n(-\f 12)=\f 2{(-2)^n}-\f{E_n}{2^n}$ and so
$E_{p-3}(-\f 12)=\f{2-E_{p-3}}{2^{p-3}}\e 8-4E_{p-3}\mod p$.
By [7, (6.3)], for $x\in\Bbb Z_{(p)}$ and $x'=(x-\xp)/p$,
\begin{align*} G_{p-1}(x)&\e (-1)^{\xp}+p^2x'(x'+1)E_{p-3}(-x)
-2p(-1)^{\xp}\f{B_{p^2(p-1)}(-x)-B_{p^2(p-1)}}{p^2(p-1)}\\&\q+2p^2(-1)^{\xp}
\Big(\Ls{B_{p^2(p-1)}(-x)-B_{p^2(p-1)}}{p^2(p-1)}^2
+\f 12(-1)^{\xp}E_{p-3}(-x)\Big)\mod{p^3}.\end{align*}
Taking $x=\f 12$ and then applying (3.2) yields that
\begin{align*}G_{p-1}\Ls 12&\e (-1)^{\f{p+1}2}-\f 14p^2E_{p-3}\Big(-\f 12\Big)-2p(-1)^{\f{p+1}2}(2-2\qp 2+p\qp 2^2)\\&\q+2p^2(-1)^{\f{p+1}2}(2-2\qp 2)^2+p^2E_{p-3}\Big(-\f 12\Big)
\\&\e (-1)^{\f{p+1}2}-4(-1)^{\f{p+1}2}(1-\qp 2)p\\&\q+p^2(6-3E_{p-3}+(-1)^{\f{p+1}2}(8-16\qp 2+6\qp 2^2))
\mod {p^3}.
\end{align*}
It is clear that
\begin{align*}\b{\f 12}p\b{-\f 32}p&=\f{1+2p}{1-2p}\b{-\f 12}p^2
=\f{1+2p}{1-2p}\cdot\f{\b{2p}p^2}{16^p}\e \f{1+2p}{1-2p}\cdot\f 1{4\cdot 16^{p-1}}
\\&\e \f{(1+2p)(1+2p+4p^2)}{4\cdot 16^{p-1}}\e \f{1+4p+8p^2}{4\cdot 16^{p-1}}\mod {p^3}
\end{align*}
and
\begin{align*}
&16^{p-1}(1+4p(1-\qp 2)+p^2(8-16\qp 2+10\qp 2^2))
\\\e& (1+4p\qp 2+6p^2\qp 2^2)(1+4p(1-\qp 2)+p^2(8-16\qp 2+10\qp 2^2))
\\\e& 1+4p+8p^2\mod {p^3}.
\end{align*}
By [7, (6.6)], for $x\in\Bbb Z_{(p)}$,
\begin{align*} G_p(x)&\e 1-\b{x}p\b{-1-x}p+
2p\f{B_{p^2(p-1)}(-x)-B_{p^2(p-1)}}{p^2(p-1)}
\\&\q +2p^2\Big(\f{B_{p^2(p-1)}(-x)-B_{p^2(p-1)}}{p^2(p-1)}\Big)^2
+p^2(-1)^{\xp}E_{p-3}(-x)\mod {p^3}.
\end{align*}
Taking $x=\f 12$ and then applying (3.2) and the above we deduce that
\begin{align*}G_p\ls 12&\e 1-\b{\f 12}p\b{-\f 32}p+2p(2-2\qp 2+p\qp 2^2)
\\&\q+2p^2(2-2\qp 2)^2+p^2(-1)^{\f{p+1}2}E_{p-3}\Big(-\f 12\Big)
\\&\e 1-\f{1+4p+8p^2}{4\cdot 16^{p-1}}
+4p(1-\qp 2)+p^2(8-16\qp 2+10\qp 2^2)\\&\q+(-1)^{\f{p+1}2}p^2(8-4E_{p-3})\mod{p^3}
\end{align*}
and so
$$16^{p-1}G_p\Ls 12\e \Big(1-\f 14\Big)(1+4p+8p^2)+(-1)^{\f{p+1}2}p^2(8-4E_{p-3})\mod {p^3}.$$
Hence,
\begin{align*}
G_p\Ls 12G_{p-1}\Ls 12&\e \Big(16^{-(p-1)}\cdot \f 34(1+4p+8p^2)+(-1)^{\f{p+1}2}p^2(8-4E_{p-3})\Big)
\\&\q\times\Big((-1)^{\f{p+1}2}-4(-1)^{\f{p+1}2}(1-\qp 2)p\\&\q+p^2(6-3E_{p-3}+(-1)^{\f{p+1}2}(8-16\qp 2+6\qp 2^2))\Big)
\\&\e\f 34(-1)^{\f{p+1}2}\f{1+4p\qp 2+6p^2\qp 2^2}{16^{p-1}}+\Big(\f{25}2-\f{25}4E_{p-3}\Big)p^2
\\&\e \f 34(-1)^{\f{p+1}2}+p^2\Big(\f{25}2-\f{25}4E_{p-3}\Big)
\mod {p^3}.\end{align*}
Since
$G_n\sls 12=g_n\sls 52$, from Corollary 2.1 and the above we derive
\begin{align*}
\sum_{n=0}^{p-1}(2n+1)G_n\Ls 12^2
&=p^2G_p\Ls 12G_{p-1}\Ls 12-\f 43p^4\Big(G_p\Ls 12-G_{p-1}\Ls 12\Big)^2
\\&\e p^2\Big(\f 34(-1)^{\f{p+1}2}+p^2\Big(\f{25}2-\f{25}4E_{p-3}\Big)\Big)
-\f 43p^4\Big(\f 34-(-1)^{\f{p+1}2}\Big)^2
\\&=\f 34(-1)^{\f{p+1}2}p^2+\Big(\f{125}{12}+2(-1)^{\f{p+1}2}-\f{25}4
E_{p-3}\Big)p^4\mod {p^5}.
\end{align*}
This completes the proof.
\vskip0.2cm
\par{\bf Theorem 3.11} {\sl Let $p>3$ be a prime, $x\in\Bbb Z_{(p)}$, $x\not\e 1,5\mod p$ and $\sls{2x-1}p=1$.
Then
\begin{align*}&\sum_{n=0}^{p-1}(2n+1)^3g_n(x)^2
\\&\e \f{x+2}{3}(-1)^{x_0}p^2
\Big(1+\f{x_1(x_1+8x_0-4)}{16(2x-1)}
+\f{x_1^2(x_1+12x_0-6)}{32(2x-1)^2}p\Big)\mod {p^4},
\end{align*}
where $x_0$ and $x_1$ are given by (3.1)}
\vskip0.2cm
\par{\it Proof.} From Theorem 2.7 and (3.3),
\begin{align*}&3(x-1)(x-5)\sum_{n=0}^{p-1}(2n+1)^3g_n(x)^2
\\&\e (x-1)(x^2-3x-10)p^2 (-1)^{x_0}\Big(1+\f{x_1(x_1+8x_0-4)}{16(2x-1)}
+\f{x_1^2(x_1+12x_0-6)}{32(2x-1)^2}p\Big)\mod {p^4}.
\end{align*}
This yields the result.
\vskip0.2cm
\par{\bf Corollary 3.1} {\sl Let $p>3$ be a prime, $x\in\Bbb Z_{(p)}$, $x\not\e -2,-1,0,1,-\f 12\mod p$ and $x'=(x-\xp)/p$. Then}
$$\sum_{n=0}^{p-1}(2n+1)^3G_n(x)^2\e \f{2x^2+2x+3}3
(-1)^{\xp-1}(1+x'(x'+1))p^2\mod {p^4}.$$
\par{\it Proof.}
Set
$$ x_0=\begin{cases}\xp+1=x+1-px'&\t{if $\xp<\f{p-1}2$,}
\\p-\xp=(x'+1)p-x&\t{if $\xp\ge \f{p-1}2$}
\end{cases}$$
and
$$ x_1=\f{(2x+1)^2-(2x_0-1)^2}p=\begin{cases}
4x'(2x+1-px')&\t{if $\xp<\f{p-1}2$,}
\\ 4(x'+1)(2x+1-p(x'+1))&\t{if $\xp\ge \f{p-1}2$.}
\end{cases}$$
Then $x_0\in\{1,2,\ldots,\f{p-1}2\}$ and
$2(2x^2+2x+1)-1=(2x+1)^2=(2x_0-1)^2+px_1$.
Since $G_n(x)=g_n(2x^2+2x+1)$, replacing $x$ with $2x^2+2x+1$ in Theorem 3.11 gives
\begin{align*}&\sum_{n=0}^{p-1}(2n+1)^3G_n(x)^2
\\&\e\f{2x^2+2x+1+2}{3}
(-1)^{x_0}p^2\Big(1+\f{x_1(x_1+8x_0-4)}{16(2x+1)^2}
+\f{x_1^2(x_1+12x_0-6)}{32(2x+1)^4}p\Big)
\mod {p^4}.
\end{align*}
By the proof of Theorem 3.5,
$$(-1)^{x_0}p^2\Big(1+\f{x_1(x_1+8x_0-4)}{16(2x+1)^2}
+\f{x_1^2(x_1+12x_0-6)}{32(2x+1)^4}p\Big)\e -(-1)^{\xp}
(1+x'(x'+1))p^2\mod {p^4}.$$
Thus the result follows.
\vskip0.2cm
\par{\bf Theorem 3.12} {\sl Let $p>11$ be a prime. Then}
\begin{align*}
&\sum_{n=0}^{p-1}(2n+1)^3G_n\Ls 12^2\e \f 98(-1)^{\f{p-1}2}p^2+\Big(-\f{1145}{72}+\f{14}{15}
(-1)^{\f{p-1}2}+\f{75}{8}E_{p-3}\Big)p^4\mod {p^5},
\\&\sum_{n=0}^{p-1}(2n+1)^3G_n\Big(-\f 12\Big)^2\e
-\f{5}{8}\s2 p^2+\Big(-\f{625}{216}+\f{38}9\s2-\f{125}{24}E_{p-3}\Big)p^4 \mod {p^5},
\\&\sum_{n=0}^{p-1}(2n+1)^3G_n\Big(-\f 13\Big)^2
\\&\q\e
-\f{161}{243}(-1)^{[\f p3]}p^2 +\Big(-\f{1703}{486}+\f{49}{10}(-1)^{[\f p3]}-\f{2990}{243}U_{p-3}\Big)p^4 \mod {p^5},
\\&\sum_{n=0}^{p-1}(2n+1)^3G_n\Big(-\f 14\Big)^2
\\&\q\e -\f{91}{128}(-1)^{[\f p4]}p^2+ \Big(-\f{36635}{8064}+\f{1898}{315}(-1)^{[\f p4]}-\f{2975}{128}s_{p-3}\Big)p^4\mod {p^5},
\\&\sum_{n=0}^{p-1}(2n+1)^3G_n\Big(-\f 16\Big)^2
\\&\q\e
-\f{1519}{1944}\s2 p^2 +\Big(-\f{5184329}{748440}+\f{3286}{385}\s2-\f{552965}{17496}
E_{p-3}\Big)p^4 \mod {p^5}.
\end{align*}
\par{\it Proof.} From Corollary 2.3, for $x\in\Bbb Z_{(p)}$,
\begin{align*}
&12x(x^2-1)(x+2)\sum_{n=0}^{p-1}(2n+1)^3G_n(x)^2
\\&\e 2p^4(4x^4+8x^3-2x^2-6x-8)(G_p(x)^2
+G_{p-1}(x)^2)
+(-8p^4(3x^2+3x-4)
\\&\q-4p^2x(x^2-1)(x+2)(2x^2+2x+3))
G_p(x)G_{p-1}(x)\mod {p^5}.\tag{3.6}
\end{align*}
Taking $x=\f 12$ in (3.6) and then applying the congruences for $G_p\sls 12,G_{p-1}\sls 12$ and $G_p\sls 12G_{p-1}\sls 12$ modulo $p^3$ yields
\begin{align*}&-\f{45}4\sum_{n=0}^{p-1}(2n+1)^3G_n\Ls 12^2
\\&\e -\f{41}2p^4\Big(G_p\Ls 12^2+G_{p-1}\Ls 12^2\Big)
+p^2\Big(14p^2+\f{135}8\Big)G_p\Ls 12G_{p-1}\Ls 12
\\&\e -\f{41}2p^4\Big(\f 9{16}+1\Big)+14p^4\cdot \f 34(-1)^{\f{p+1}2}
+\f{135}8p^2
\Big(\f 34(-1)^{\f{p+1}2}+p^2\Big(\f{25}2-\f{25}4E_{p-3}\Big)\Big)
\\&=\f{405}{32}(-1)^{\f{p+1}2}p^2+\Big(\f{5725}{32}+\f{21}2
(-1)^{\f{p+1}2}-\f{3375}{32}E_{p-3}\Big)p^4\mod {p^5}
\end{align*}
and so
$$\sum_{n=0}^{p-1}(2n+1)^3G_n\Ls 12^2\e
\f 98(-1)^{\f{p-1}2}p^2+\Big(-\f{1145}{72}+\f{14}{15}
(-1)^{\f{p-1}2}+\f{75}{8}E_{p-3}\Big)p^4 \mod {p^5}.$$
Taking $x=-\f 12$ in (3.6) and then applying Lemmas 3.6, 3.7 and 3.9 yields
\begin{align*}&\f{27}4\sum_{n=0}^{p-1}(2n+1)^3G_n\Big(-\f 12\Big)^2
\\&\e -\f{25}2p^4\Big(G_p\Big(-\f 12\Big)^2+G_{p-1}\Big(-\f 12\Big)^2\Big)
+p^2\Big(38p^2-\f{45}8\Big)G_p\Big(-\f 12\Big)G_{p-1}\Big(-\f 12\Big)
\\&\e -\f{25}2p^4\Big(\f 9{16}+1\Big)+38p^4\cdot \f 34(-1)^{\f{p-1}2}
-\f{45}8p^2\Big(\f 34(-1)^{\f{p-1}2}+\f{25}4p^2E_{p-3}\Big)
\\&=-\f{135}{32}\s2 p^2+\Big(-\f{625}{32}+\f{57}2\s2-\f{1125}{32}E_{p-3}\Big)p^4\mod {p^5}
\end{align*}
and so
$$\sum_{n=0}^{p-1}(2n+1)^3G_n\Big(-\f 12\Big)^2
\e -\f{5}{8}\s2 p^2+\Big(-\f{625}{216}+\f{38}9\s2-\f{125}{24}E_{p-3}\Big)p^4\mod {p^5}.$$
Taking $x=-\f 13$ in (3.6) and then applying Lemmas 3.6, 3.8 and 3.9 yields
\begin{align*}&\f{160}{27}\sum_{n=0}^{p-1}(2n+1)^3G_n\Big(-\f 13\Big)^2
\\&\e -\f{1048}{81}p^4\Big(G_p\Big(-\f 13\Big)^2+G_{p-1}\Big(-\f 13\Big)^2\Big)
+p^2\Big(\f{112}3p^2-\f{3680}{729}\Big)G_p\Big(-\f 13\Big)G_{p-1}\Big(-\f 13\Big)
\\&\e -\f{1048}{81}p^4\Big(\f {49}{81}+1\Big)+\f{112}3p^4\cdot \f 79(-1)^{[\f p3]}
-\f{3680}{729}p^2\Big(\f 79(-1)^{[\f p3]}+\f{130}9p^2U_{p-3}\Big)
\\&=-\f{25760}{6561}(-1)^{[\f p3]}p^2 +\Big(-\f{136240}{6561}+\f{784}{27}(-1)^{[\f p3]}-\f{478400}{6561}U_{p-3}\Big)p^4\mod {p^5}
\end{align*}
and so
$$\sum_{n=0}^{p-1}(2n+1)^3G_n\Big(-\f 13\Big)^2\e
-\f{161}{243}(-1)^{[\f p3]}p^2 +\Big(-\f{1703}{486}+\f{49}{10}(-1)^{[\f p3]}-\f{2990}{243}U_{p-3}\Big)p^4\mod {p^5}.$$
Taking $x=-\f 14$ in (3.6) and then applying Lemmas 3.6, 3.8 and 3.9 yields
\begin{align*}&\f{315}{64}\sum_{n=0}^{p-1}(2n+1)^3G_n\Big(-\f 14\Big)^2
\\&\e -\f{431}{32}p^4\Big(G_p\Big(-\f 14\Big)^2+G_{p-1}\Big(-\f 14\Big)^2\Big)
+p^2\Big(\f{73}2p^2-\f{2205}{512}\Big)G_p\Big(-\f 14\Big)G_{p-1}\Big(-\f 14\Big)
\\&\e -\f{431}{32}p^4\Big(\f {169}{256}+1\Big)+\f{73}2p^4\cdot \f {13}{16}(-1)^{[\f p4]}
-\f{2205}{512}p^2\Big(\f {13}{16}(-1)^{[\f p4]}+\f{425}{16}p^2s_{p-3}\Big)
\\&=-\f{28665}{8192}(-1)^{[\f p4]}p^2 +\Big(-\f{183175}{128}+1898(-1)^{[\f p4]}-\f{937125}{128}s_{p-3}\Big)\f{p^4}{64}\mod {p^5}
\end{align*}
and so
\begin{align*}&\sum_{n=0}^{p-1}(2n+1)^3G_n\Big(-\f 14\Big)^2
\\&\e -\f{91}{128}(-1)^{[\f p4]}p^2+ \Big(-\f{36635}{8064}+\f{1898}{315}(-1)^{[\f p4]}-\f{2975}{128}s_{p-3}\Big)p^4
\mod {p^5}.\end{align*}
Finally, putting $x=-\f 16$ in (3.6) and then applying Lemmas 3.6, 3.8 and 3.9 yields
\begin{align*}&\f{385}{108}\sum_{n=0}^{p-1}(2n+1)^3G_n\Big(-\f 16\Big)^2
\\&\e -\f{2297}{162}p^4\Big(G_p\Big(-\f 16\Big)^2+G_{p-1}\Big(-\f 16\Big)^2\Big)
+p^2\Big(\f{106}3p^2-\f{18865}{5832}\Big)G_p\Big(-\f 16\Big)G_{p-1}\Big(-\f 16\Big)
\\&\e -\f{2297}{162}p^4\Big(\f {961}{1296}+1\Big)+\f{106}3p^4\cdot \f {31}{36}\s2
-\f{18865}{5832}p^2\Big(\f {31}{36}\s2+\f{11285}{324}p^2E_{p-3}\Big)
\\&=-\f{385\cdot 1519}{108\cdot 1944}\s2 p^2 +\Big(-\f{5184329}{748440}+\f{3286}{385}\s2-\f{552965}{17496}
E_{p-3}\Big)\f{385p^4}{108}\mod {p^5}
\end{align*}
and so
\begin{align*}&\sum_{n=0}^{p-1}(2n+1)^3G_n\Big(-\f 16\Big)^2
\\&\e -\f{1519}{1944}\s2 p^2 +\Big(-\f{5184329}{748440}+\f{3286}{385}\s2-\f{552965}{17496}
E_{p-3}\Big)p^4\mod {p^5}.\end{align*}
This completes the proof.
\vskip0.2cm
\par{\bf Theorem 3.13} {\sl Let $p$ be an odd prime, $x\in\Bbb Z_{(p)}$, $x\not\e \f 52,\f {17}2\mod p$ and $\ls{2x-1}p=1$. Then
\begin{align*}&\sum_{n=0}^{p-1}(2n+1)^4g_n(x)^2
\\&\e
(-1)^{x_0-1}\f{3x^4-20x^3-78x^2+268x-5}{8(2x-5)(2x-17)}\Big(\f{x_1+4x_0-2}{2(2x-1)}p
+\f{x_1(x_1+8x_0-4)}{8(2x-1)^2}p^2\Big)
\mod {p^3},\end{align*}
where $x_0$ and $x_1$ are given by (3.1). }
\vskip0.2cm
\par{\it Proof.} By Theorem 2.8,
$$\sum_{n=0}^{p-1}(2n+1)^4g_n(x)^2\e \f{3x^4-20x^3-78x^2+268x-5}{8(2x-5)(2x-17)}\sum_{n=0}^{p-1}g_n(x)^2\mod {p^3}.\eqno{(3.7)}$$
This together with Lemma 3.5 yields the result.
\vskip0.2cm
\par{\bf Theorem 3.14} {\sl Let $p$ be an odd prime, $x\in\Bbb Z_{(p)}$, $x\not\e -\f 52,\pm\f 32,\pm \f 12\mod p$ and $x'=(x-\xp)/p$. Then
\begin{align*}&\sum_{n=0}^{p-1}(2n+1)^4G_n(x)^2
\\&\e (-1)^{\xp}\f{(6x^8+24x^7+28x^6-78x^4-128x^3-44x^2+16x+21)
(1+2x')}{(4x^2-1)(4x^2-9)(2x+5)}p
\mod {p^3}.\end{align*}}
\par{\it Proof.} Since $G_n(x)=g_n(2x^2+2x+1)$, the result follows from (3.7) and (1.2).
\vskip0.2cm

\par {\bf Theorem 3.15} {\sl Let $p$ be a prime greater than  $3$. Then}
$$\sum_{n=0}^{p-1}G_n\Ls 32^2=
\sum_{n=0}^{p-1}G_n\Big(-\f 52\Big)^2=\sum_{n=0}^{p-1}
g_n\Ls{17}2^2\e -\f{10}3\s2p^3+\f{91}{27}p^4\mod {p^5}.$$
\par{\it Proof.} Taking $x=\f{17}2$ in Theorem 2.8 gives
$$-\f{243}{16}\sum_{n=0}^{p-1}g_n\Ls{17}2^2
\e 117p^4\Big(g_p\Ls{17}2^2-g_{p-1}\Ls{17}2^2\Big)+\f{135}2p^3
g_p\Ls{17}2g_{p-1}\Ls{17}2\mod {p^5}.$$
Since $2\cdot \f{17}2-1=(2\cdot\f{p-3}2-1)^2+p(8-p)$,
for $x=\f{17}2$ we have $x_0=\f{p-3}2$ and $x_1=8-p$ by (3.1).
Hence, from Lemma 3.4 and (3.3) we see that
$$g_p\Ls{17}2\e 1+\f{(8-p)(8-p+4(p-3)-4)}{16^2}
\e \f 34\mod p,\  g_{p-1}\Ls{17}2\e (-1)^{\f{p-1}2}\mod p$$
and
\begin{align*}&g_{p}\Ls{17}2g_{p-1}\Ls{17}2
\\&\e \s2\Big(1+\f{(8-p)(8-p+4(p-3)-4)}{16^2}+\f{(8-p)^2(8-p+6(p-3)-6)}
{32(17-1)^2}p\Big)
\\&\e \f 34(-1)^{\f{p-1}2}\mod {p^2}.
\end{align*}
Thus,
\begin{align*}\sum_{n=0}^{p-1}g_n\Ls{17}2^2
&\e -\f{16}{243}\Big(117p^4\Big(\f 9{16}-1\Big)+\f{135}2p^3
\cdot\f 34\s2\Big)
\\&=-\f{10}3\s2 p^3+\f{91}{27}p^4\mod {p^5}.
\end{align*}
To complete the proof, we note that $G_n\sls 32=G_n(-\f 52)=g_n\sls{17}2$ since $G_n(x)=g_n(2x^2+2x+1)$.


\begin{thebibliography}{99}
 \bibitem [1] {} V.J.W. Guo, {\it Proof of Sun's conjectures on
integer-valued polynomials}, J. Math. Anal. Appl. {\bf 444} (2016), 182-191.

       \bibitem [2] {} J.-C. Liu, {\it Proof of some divisibility results on sums involving binomial coefficients,} J. Number Theory {\bf 180} (2017), 566-572.

         \bibitem [3] {} J.-C. Liu, {\it A generalized supercongruence of Kimoto and Wakayama}, J. Math. Anal. Appl. {\bf 467} (2018), 15-25.

              \bibitem [4] {} J.-C. Liu and H.-X. Ni, {\it On two supercongruences involving
Almkvist-Zudilin sequences}, Czech. Math. J. {\bf 71}(2021), 1211-1219.

\bibitem [5] {} L. Long, R. Osburn and H. Swisher, {\it On a conjecture of Kimoto and Wakayama}, Proc. Amer. Math. Soc.
    {\bf 144} (2016), 4319-4327.

\bibitem [6] {} W. Magnus, F. Oberhettinger and R. P. Soni,
{\it Formulas and Theorems for the Special Functions of Mathematical Physics}, 3rd edn.
Springer, New York, 1966.

\bibitem [7] {} Z.H. Sun, {\it
Congruences for certain families of Ap\'ery-like sequences}, Czech.
Math. J. {\bf 72} (2022), 875-912.

\bibitem [8] {}  Z.H. Sun, {\it Binomial Coefficients, Recurrence Sequences and Congruences} (Chinese), Science Press, Beijing, 2025.

\bibitem [9] {} Z. H. Sun, {\it Congruences for a type of Ap\'ery-like numbers}, Chin. Ann. Math., Ser. B, to appear.

    \bibitem [10] {} Z. H. Sun, {\it Generalizations of the Christoffel-Darboux formula  and congruences involving Ap\'ery-like numbers}, arXiv:2608.13192.

    \bibitem [11] {} Z.W. Sun, {\it
Supercongruences involving dual sequences}, Finite Fields Appl. {\bf
46} (2017), 179-216.

      \bibitem [12] {} C. Wang and S.-J. Wang, {\it On a conjectural supercongruence involving the dual sequence $s_n(x)$}, Proc. Amer. Math. Soc. {\bf 154} (2026), 3635-3649.

\end{thebibliography}
\end{document}